\documentclass[12pt]{article}
\usepackage{amssymb}
\usepackage{latexsym}
\def\part#1{\frac{\partial\phantom{q}}{\partial#1}}

\newenvironment{rmk}{\begin{trivlist}\item[]{\bf Remark:} }
{\end{trivlist}}
\newenvironment{ex}{\begin{trivlist}\item[]{\bf Example:} }
{\end{trivlist}}
\newenvironment{prf}{\begin{trivlist}\item[]{\bf Proof:} }
{\hfill $\Box$ \end{trivlist}}

\newtheorem{thm}{Theorem}
\newtheorem{definition}{Definition}
\newtheorem{prp}[thm]{Proposition}

\newcommand{\lie}[1]{\mathfrak{#1}}
\def\End{\mathop{\rm End}\nolimits}

\def\tr{\mathop{\rm tr}\nolimits}

\def\Diff{\mathop{\rm Diff}\nolimits}

\newcommand{\R}{\mathbf{R}}
\newcommand{\C}{\mathbf{C}}

\newcommand{\Z}{\mathbf{Z}}

\begin{document}
\title{Brackets, forms and invariant functionals}
 \author{Nigel Hitchin\\[5pt]}
\maketitle

\centerline{\it Dedicated to the memory of Shiing-Shen Chern}

\begin{abstract}
\noindent In the context of generalized geometry we first show how the Courant bracket helps to define connections with skew torsion and then investigate a five-dimensional invariant functional and its associated geometry, which involves three Courant-commuting sections of $T\oplus T^*$. A Hamiltonian flow arising from this corresponds to a version of the Nahm equations, and we investigate the six-dimensional geometrical structure this describes.
 \end{abstract}
\section{Introduction}
Chern was an acknowledged master in the use of exterior differential forms to reveal  geometrical truths. In this paper we shall also use differential forms,  adopting however a somewhat unorthodox point of view, initiated in the author's paper \cite{Hit3}. We regard forms as {\it spinors} for the bundle $T\oplus T^*$ -- a section $X+\xi$ of this bundle acts on a form $\rho$ by $(X+\xi)\cdot \rho=i_X\rho+\xi\wedge \rho$ and this satisfies the Clifford algebra identities for the indefinite metric on $T\oplus T^*$ defined by $(X+\xi,X+\xi)=i_X\xi$. Once we look at forms this way, then a natural operation on sections of $T\oplus T^*$ appears -- the {\it Courant bracket}. 

In the first part of this paper we show how the Courant bracket can be used to define some familiar objects in differential geometry -- in particular the Levi-Civita connection and more generally connections with skew torsion, a topic which has had prominence recently because of its use in String Theory.

The second part uses the invariant functional approach of \cite{Hit3} to define and study  geometrical structures in five and six dimensions which are defined  by closed differential forms. Their characteristic feature consists of a  triple of sections of $T\oplus T^*$. In the five-dimensional case this is a Courant-commuting triple, and in six dimensions a solution to an analogue of Nahm's equations with the Courant bracket replacing the bracket of a finite-dimensional Lie algebra. What we obtain is a novel six-dimensional structure which has some of the features of four-dimensional self-duality.

\section{Generalized geometry}

\subsection{The basic setting}
The essential idea is to take a manifold $M^n$ and replace the tangent bundle $T$ by $T\oplus T^*$. This has a natural inner product of signature $(n,n)$ defined by
$$(X+\xi,X+\xi)=i_X\xi$$
for a tangent vector $X$ and cotangent vector $\xi$ (we adhere to the conventional choice of sign here rather than its negative as in \cite{Hit3}). The skew-adjoint endomorphisms of $T\oplus T^*$ are sections of the bundle $\End T\oplus \Lambda^2T^*\oplus \Lambda^2T$ and we focus in particular on the action of a $2$-form $B$. Exponentiating it to an orthogonal transformation of $T\oplus T^*$, we get $$X+\xi\mapsto X+\xi+i_XB.$$
The bundle of differential forms $\Lambda^{\ast}T^*$ we consider as a bundle of Clifford modules over the Clifford algebra generated by the action of $T\oplus T^*$:
$$(X+\xi)\cdot\varphi=i_X\varphi+\xi\wedge \varphi.$$
This satisfies the relation $(X+\xi)^2=(X+\xi,X+\xi)1$. The $2$-form $B$ when exponentiated into the spin group acts on a form as 
$$\varphi\mapsto e^B\varphi.$$
In a sense, the exterior derivative is dual to the Lie bracket -- one can be defined in terms of the other and which comes first in an elementary  course depends on the lecturer's preference for forms or vector fields. Here is one version of the relation:
$$2i_{[X,Y]}\alpha=d([i_X,i_Y]\alpha)+2i_Xd(i_Y\alpha)-2i_Yd(i_X\alpha)+[i_X,i_Y]d\alpha$$
In our current view we treat forms as spinors and then replacing the vector fields $X,Y$ in this formula by sections $u,v$ of $T\oplus T^*$, and interior products by Clifford products, we  define a  bracket $[u,v]$ by
\begin{equation}
2[u,v]\cdot\alpha=d((u\cdot v-v \cdot u)\cdot\alpha)+2u\cdot d(v\cdot\alpha)-2v\cdot d(u\cdot\alpha)+(u\cdot v-v \cdot u)\cdot  d\alpha
\label{courant}
\end{equation}
This is the {\it Courant bracket} \cite{C} defined explicitly as 
\begin{equation}
[X+\xi, Y+\eta]=[X,Y]+\mathcal{L}_{X}\eta -\mathcal{L}_{Y}\xi -\frac{1}{2}
  d(i_{X}\eta  -i_{Y}\xi).
  \label{cdef}
  \end{equation}
  It satisfies the two basic identities
  \begin{eqnarray}
  [u,fv]&=&f[u,v]+(\pi(u)f)v-(u,v)df \label{c1}\\
  \pi(u)(v,w)&=&([u,v]+d(u,v),w)+(v,[u,w]+d(u,w))\label{c2}
  \end{eqnarray}
  where $\pi(X+\xi)=X$.
  
  If $B$ is a {\it closed} $2$-form, then $e^B$ commutes with $d$ and so, by  its definition using the exterior derivative, the Courant bracket is invariant under the map $X+\xi\mapsto X+\xi+i_XB$. 

``Generalized geometries" consist of data on $T\oplus T^*$ which are compatible with the pre-existing $SO(n,n)$ structure and satisfy an integrability condition which is expressed using  the exterior derivative or the Courant bracket. They can be transformed by an element of the semi-direct product $\Diff(M) \ltimes \Omega_{closed}^2(M)$. This group replaces the diffeomorphism group. Note from (\ref{cdef}) that $X+\xi\mapsto X+d\xi$ takes the Courant bracket to the bracket in the Lie algebra of this group.

\begin{ex} To incorporate a Riemannian metric $g$ into this picture we consider it as a homomorphism $g:T\rightarrow T^*$ and take its graph in $T\oplus T^*$. This is a subbundle $V\subset T\oplus T^*$ on which the indefinite inner product is positive definite.
\end{ex}

\subsection{Gerbes}

The rationale for treating forms as spinors, forgetting the $\Z$-grading and the algebra structure, may be considered artificial but it is more natural when we twist the picture with a {\it gerbe}. Here we adopt the naive point of view of gerbes advanced in \cite{Hit4}.

We take an open covering $\{U_{\alpha}\}$ of a manifold $M$ and think of a gerbe on $M$ in terms of a 2-cocycle: that is a collection of functions 
 $$g_{\alpha\beta\gamma}:U_\alpha\cap U_\beta\cap U_\gamma\rightarrow S^1$$
 where $g_{\alpha\beta\gamma}=g^{-1}_{\beta\alpha\gamma}=\dots$ and 
 $\delta g=g_{\beta\gamma\delta}g^{-1}_{\alpha\gamma\delta}g_{\alpha\beta\delta}g^{-1}_{\alpha\beta\gamma}=1$ 
   on $U_\alpha\cap U_\beta\cap U_\gamma \cap U_\delta$.
 \vskip .25cm
A {\it trivialization} of a gerbe is given by functions $f_{\alpha\beta}: U_\alpha\cap U_\beta\rightarrow S^1$ with $f_{\alpha\beta}=f^{-1}_{\beta\alpha}$ and satisfying $$g_{\alpha\beta\gamma}=f_{\alpha\beta}f_{\beta\gamma}f_{\gamma\alpha}.$$ 
Take two trivializations, then
$g_{\alpha\beta\gamma}=f_{\alpha\beta}f_{\beta\gamma}f_{\gamma\alpha}=\tilde f_{\alpha\beta}\tilde f_{\beta\gamma}\tilde f_{\gamma\alpha}$,
and so 
the ratio of two trivializations $h_{\alpha\beta}=(f/\tilde f)_{\alpha\beta}$ defines the transition functions for  a principal circle bundle:
$$h_{\alpha\beta}h_{\beta\gamma}h_{\gamma\alpha}=1.$$
Following \cite{Bry}, we have 
\begin{definition}
A {\rm connective structure} on a gerbe is defined by a collection of $1$-forms $A_{\alpha\beta}\in \Omega^1(U_\alpha\cap U_\beta)$ where
$A_{\alpha\beta}+A_{\beta\gamma}+A_{\gamma\alpha}
= g_{\alpha\beta\gamma}^{-1}dg_{\alpha\beta\gamma}$
on $U_{\alpha}\cap U_{\beta}\cap U_{\gamma}$.
\end{definition}

\begin{definition}
A {\rm curving} of a connective structure  is defined by a collection of $2$-forms $B_{\alpha}\in \Omega^2(U_\alpha)$ where
$B_\beta-B_\alpha=dA_{\alpha\beta}$ on $U_{\alpha}\cap U_{\beta}$. Then
 $dB_{\beta}=dB_\alpha$ is a global closed three-form $H$, called the curvature.
 \end{definition}
\vskip .25cm
\begin{rmk} A trivialization $h_{\alpha\beta}$  is {\it flat } for a connective structure if
$$h^{-1}_{\alpha\beta}dh_{\alpha\beta}=A_{\alpha\beta}.$$
The ratio of two flat trivializations is a line bundle with constant transition functions -- a flat line bundle. Now given a map $f:S^1\rightarrow M$, we pull back the gerbe to the circle, where any gerbe has a flat trivialization. Let us identify two flat trivializations on the circle if their ratio is a flat $S^1$ bundle with trivial holonomy. Then the equivalence classes define  a space on which the unit circle acts transitively. A connective structure on a gerbe then defines a principal circle bundle on the loop space of $M$: a curving defines a connection on it whose curvature transgresses $H$.
\end{rmk}

\subsection{The generalized tangent bundle}

Given a connective structure, we have
$$A_{\alpha\beta}+A_{\beta\gamma}+A_{\gamma\alpha}
= g_{\alpha\beta\gamma}^{-1}dg_{\alpha\beta\gamma}$$
and hence
$$dA_{\alpha\beta}+dA_{\beta\gamma}+dA_{\gamma\alpha}
= 0.$$
This  is a cocycle with values in closed $2$-forms and defines a vector bundle as an  extension
$$0\rightarrow T^*\rightarrow E\rightarrow T\rightarrow 0$$
obtained by identifying $T\oplus T^*$ on $U_{\alpha}$ with $T\oplus T^*$ on $U_{\beta}$ by the B-field action $X+\xi\mapsto X+\xi+i_XdA_{\alpha\beta}$.

The action of the two-form $dA_{\alpha\beta}$ on $T\oplus T^*$   preserves the natural $SO(n,n)$ structure, so $E$ also has such a structure, and since $dA_{\alpha\beta}$ is closed, it preserves  the Courant bracket, 
  so there is an induced bracket on sections of $E$ satisfying (\ref{c1}), (\ref{c2}), where $\pi:E\rightarrow T$ is the projection. We call this the {\it generalized tangent bundle}.
  
 If we now look at the action on forms $\varphi\mapsto e^{dA_{\alpha\beta}}\varphi$, then  we obtain another bundle -- the spinor bundle for $E$, which still admits the exterior derivative $d$, since this commutes with $e^{dA_{\alpha\beta}}$, but it has lost its $\Z$-grading. The cohomology of $d:C^{\infty}(S)\rightarrow C^{\infty}(S)$ is the well-known {\it twisted cohomology}. A  section of $S$ is represented by forms $\varphi_{\alpha}$ such that on $U_{\alpha}\cap U_{\beta}$ we have 
$$\varphi_{\alpha}=e^{dA_{\alpha\beta}}\varphi_{\beta}.$$
Given a curving, we have $B_\beta-B_\alpha=dA_{\alpha\beta}$ and so
$$e^{B_{\alpha}}\varphi_{\alpha}=e^{B_{\beta}}\varphi_{\beta}.$$
This defines a global form  $\psi=e^{B_{\alpha}}\varphi_{\alpha}$ but now if $d\varphi_{\alpha}=0$,
$$(d-H)\psi=0$$
where $H=dB_{\alpha}$, which is the usual definition of twisted cohomology.

We have a choice here: with a connective structure we need to adopt the spinor viewpoint, but if we insist on reverting to forms we need the curving. Note that a $2$-form $B$ still acts on $E$: $u\mapsto u+i_{\pi(u)}B$ and on $S$. In fact the Clifford action of $T^*\subset E$ on $S$ is  exterior multiplication, so $S$ is a module over the exterior algebra.

\subsection{Connections with skew torsion}

We now generalize our description of a Riemannian metric in this set-up:

\begin{definition} Let $E$ be the generalized tangent bundle for a gerbe with connective structure. A {\rm generalized metric} is a subbundle $V\subset E$ of rank $n$ on which the induced metric is positive definite.
\end{definition}

A generalized metric defines a number of structures. Firstly, since the inner product on $V$ is positive definite, and that on $T^*$ is zero, we have $V\cap T^*=0$ and so $V$ defines a splitting of the exact sequence 
$$0\rightarrow T^*\rightarrow E\rightarrow T\rightarrow 0.$$
The inner product on $E$ has signature $(n,n)$, therefore the inner product on $V^{\perp}$ is negative definite, so this also defines a splitting. Splittings form an affine space so we can also take the average of two splittings to obtain a third. The difference of two splittings is  a homomorphism from  $T$ to $T^*$. 

\begin{prp} \label {geng} Let $V\subset E$ be a generalized metric, then
\begin{itemize}
\item 
half the  difference of the two splittings is a metric $g$ on $M$.
\item 
the average of the two splittings is a curving of the connective structure.
\end{itemize}
\end{prp}

\begin{prf} In local terms a splitting is a family of sections $C_{\alpha}$ of $T^*\otimes T^*$ which satisfy 
\begin{equation}
C_{\beta}-C_{\alpha}=dA_{\alpha\beta}
\label{split}
\end{equation}
on $U_{\alpha}\cap U_{\beta}$. Let $C_{\alpha}$ define the splitting $V$, then a tangent vector $X$ is lifted  locally to $X+i_XC_{\alpha}$ (where $i_XC(Y)$ is defined to be $C(X,Y)$. By definition, the induced inner product is positive definite so that  
$$(X+i_XC_{\alpha},X+i_XC_{\alpha})=C_{\alpha}(X,X)$$
is positive definite. 

Now $Y-i_YC^T_{\alpha}$ is orthogonal to $V$ because
$$(Y-i_YC_{\alpha}^T,X+i_XC_{\alpha})= \frac{1}{2}\left( C_{\alpha}(X,Y)-C^T_{\alpha}(Y,X)\right)=0.$$

Thus half the difference of the two splittings is  
$$\frac{1}{2}(X+i_XC_{\alpha}-X+i_XC^T_{\alpha})$$
which is the symmetric part of $C_{\alpha}$. This is positive definite and so defines a metric $g$. Another way to define it is to lift tangent vectors $X,Y$ to vectors $X^+,Y^+$ in $V$ and use the inner product: $g(X,Y)=(X^+,Y^+)$.

The average of the two splittings is the skew symmetric part of $C_{\alpha}$ --  a $2$-form $B_\alpha$ --  and from (\ref{split}) this is a curving of the connective structure. In fact it is clear from what we have done here that a curving is simply a splitting which is {\it isotropic}.
\end{prf}

Given a vector field $X$ we can lift it to a section of $E$ in two ways: $X^+$ a section of $V$ and $X^-$ a section of $V^{\perp}$. We then have the following

\begin{thm} Let $X,Y$ be two vector fields, and $V\subset E$ a generalized metric. Let $g$ be the metric on $M$ and $H$ the curvature of the gerbe as defined in Proposition \ref{geng}. Then, using the Courant bracket on sections of $E$,
$$[X^-,Y^+]-[X,Y]^-=2g\nabla_XY\in \Omega^1$$
where $\nabla$ is a connection which preserves the metric $g$  and has  skew torsion $-H$.
\end{thm}

\begin{rmk} Interchanging the roles of $V$ and $V^{\perp}$, we get a connection with torsion $+H$.
\end{rmk}
 
\begin{prf} Since $\pi[A,B]=[\pi A, \pi B]$, it is clear that $\pi([X^-,Y^+]-[X,Y]^-)=0$ and so $[X^-,Y^+]-[X,Y]^-$ is a one-form.

Put $\Delta_XY=[X^-,Y^+]-[X,Y]^-$. Then
\begin{eqnarray*}
\Delta_{fX}Y&=&[fX^-,Y^+]-[fX,Y]^-\\
&=&f[X^-,Y^+]-(Yf)X^--f[X,Y]^-+(Yf)X^-\\
&=&f\Delta_XY
\end{eqnarray*}
where we have used (\ref{c1}) together with $(X^-,Y^+)=0$.

Similarly
\begin{eqnarray*}
\Delta_{X}fY&=&[X^-,fY^+]-[X,fY]^-\\
&=&(Xf)Y^++f[X^-,Y^+] -(Xf)Y^--f[X,Y]^-\\
&=&f\Delta_XY+Xf(Y^+-Y^-)\\
&=&f\Delta_XY +(Xf)2gY
\end{eqnarray*}
where we have used half the difference of the two splittings to define the metric as in Proposition \ref{geng}.

These two expressions show that $\Delta_XY=2g\nabla_XY$ for some connection $\nabla$ on $T$. We next show it preserves the metric. 

If $\xi$ is a $1$-form then $i_Z\xi=2(\xi,Z^+)$. We therefore have
$$g(\nabla_XY,Z)+g(Y,\nabla_XZ)=(\Delta_XY,Z^+)+(Y^+,\Delta_XZ).$$
But  
$(\Delta_XY,Z^+)=([X^-,Y^+]-[X,Y]^-,Z^+)=([X^-,Y^+],Z^+)$ and so
$$(\Delta_XY,Z^+)+(Y^+,\Delta_XZ)=([X^-,Y^+],Z^+)+(Y^+,[X^-,Z^+]).$$ 
But using (\ref{c2}) on the right hand side together with  $(X^-,Y^+)=0=(X^-,Z^+)$, we find
$$(\Delta_XY,Z^+)+(Y^+,\Delta_XZ)=X(Y^+,Z^+)=Xg(Y,Z).$$

For the torsion, we consider
\begin{equation}
2g(\nabla_XY-\nabla_YX-[X,Y])=\Delta_XY-\Delta_YX-2g[X,Y]
\label{tor}
\end{equation}
Now
\begin{equation}
\Delta_XY-\Delta_YX=[X^-,Y^+]-[Y^-,X^+]-2[X,Y]^-.
\label{del}
\end{equation}

One-forms Courant-commute so 
\begin{equation}
[X^+-X^-,Y^+-Y^-]=0
\label{comm}
\end{equation} 

Consider now the Courant bracket of two vector fields lifted by the splitting defined by the curving of the gerbe. The $2$-form $B_{\alpha}$ is not in general closed, and instead of preserving the Courant bracket we get  
$$[X+i_XB_{\alpha}, Y+i_YB_{\alpha}]=[X,Y]+i_{[X,Y]}B_{\alpha}-i_Xi_YdB_{\alpha}.$$
But the curving was the average of the two splittings so this means
$$\frac{1}{4}[X^++X^-,Y^++Y^-]=\frac{1}{2}([X,Y]^++[X,Y]^-)-i_Xi_YH$$
Using (\ref{comm}) this gives 
$$[X^+,Y^-]+[X^-,Y^+]=[X,Y]^++[X,Y]^--2i_Xi_YH.$$
Thus in (\ref{del}), $\Delta_XY-\Delta_YX=[X,Y]^+-[X,Y]^--2i_Xi_YH=2g([X,Y])-2i_Xi_YH$
and from (\ref{tor}) the torsion is $-H$.
\end{prf}

\begin{ex} If we take $V$ to be the graph in $T\oplus T^*$ of a Riemannian metric, then the theorem provides,  on expanding the Courant bracket,  the familiar formula (using the summation convention) for the Levi-Civita connection:

$$\left[\frac{\partial}{\partial x_i}-g_{ik}dx_k, \frac{\partial}{\partial x_j}+g_{jk}dx_k\right]-\left[\frac{\partial}{\partial x_i},\frac{\partial}{\partial x_j}\right]^+ 
=\left(\frac{\partial g_{jk}}{\partial x_i}+\frac{\partial g_{ik}}{\partial x_j}-\frac{\partial g_{ij}}{\partial x_k}\right)dx_k=2g_{\ell k}\Gamma^{\ell}_{ij}dx_k$$
\end{ex}

\section{Invariant functionals}
\subsection{Generalized geometry from open orbits}\label{open}
The algebraic origins  for this part of the paper are laid out in the theory of {\it prehomogeneous vector spaces} as described by Kimura and Sato in \cite{Sat} or \cite{Kim}. They look for the following data:
 a Lie group $G$;
	a representation space $V$;
	an open orbit $U\subset V$, and
a relatively invariant polynomial $f$ (i.e. $f(gv)=\chi(g)f(v)$). 
				They  give a list of such irreducible representations over $\C$. Rather remarkably, the various real forms provide many situations familiar to a differential geometer, not least the lists of irreducible holonomy groups of Berger and Merkulov-Schwachh\"ofer. 

In previous papers (\cite{Hit1},\cite{Hit2}) the author has considered the variational origins of  geometrical structures derived from the members of the list which consist of $GL(n,\R)$ acting on $p$-forms. 
 In the current context, we are considering generalized geometries, which relate to open orbits of the groups $Spin(n,n)\times \R^*$ acting on spinors, which will be non-homogeneous differential forms in our realization. The procedure is as follows.

Let  $M$  be a compact oriented $n$-manifold with a closed form $\rho$ which lies in the open orbit at each point -- following \cite{Hit2} we shall call such forms {\it stable}. The invariant polynomial $f$ of degree $d$ defines $f(\rho)$ as a non-zero section of $(\Lambda^nT^*)^{d/2}$ and then $\phi(\rho)=\vert f\vert^{2/d}$ is an invariantly defined volume form. By integration we define a total volume $V(\rho)$. By openness, any nearby form will again be stable, and we look for critical points of the volume functional on a given cohomology class in $H^{\ast}(M,\R)$:
$$V(\rho)=\int_M \!\!\!\phi(\rho).$$
Now  the derivative of $\phi$ at $\rho$ is a linear map
$$D\phi:\Lambda^{\ast}T^*\rightarrow \Lambda^nT^*.$$
We can write this in terms of a form $\hat\rho$ 
$$D\phi(\dot\rho)=\langle \hat\rho,\dot\rho\rangle.$$
using the natural invariant pairing on spinors. This  is realized on forms as follows: if $\sigma (\alpha)=(-1)^m\alpha$ if $\alpha$ is of degree $2m$ or $2m+1$, then the pairing, with values in $\Lambda^nT^*$, is 
$$\langle \varphi_1,\varphi_2\rangle=[\varphi_1\wedge \sigma(\varphi_2)]_n.$$
It is often called the Mukai pairing.

A critical point of the volume functional $V$ is given by 
$$0=\delta V(\dot\rho)=\int_M\!\!\! D\phi(\dot\rho)=\int_M \!\!\langle \hat\rho,\dot\rho\rangle$$
for all exact $\dot\rho=d\sigma$. From Stokes' theorem it follows that $\rho$ must satisfy $d\hat\rho=0$. The geometry comes from  understanding the  interpretation of the two equations
\begin{equation}
d\rho=0,\qquad d\hat\rho=0.
\label{max}
\end{equation}
\vskip .25cm

The case of $ Spin(6,6)\times \R^*$, which has an open orbit with stabilizer $SU(3,3)$, gave rise to the definition of {\it generalized Calabi-Yau structure} \cite{Hit3}. Since $SU(3,3)$ preserves a complex structure, we obtain $J:T\oplus T^*\rightarrow T\oplus T^*$ with $J^2=-1$ and the integrability condition is obtained by replacing the Lie brackets in the Nijenhuis tensor by {\it Courant} brackets.

A less restrictive geometry, a {\it generalized complex structure}, is  developed in \cite{Gu}. This relies on a reduction to $U(3,3)$ with the same integrability condition.

These structures exist in all dimensions but the specific open orbit picture reveals a number of special properties such as the existence of flat coordinates on the moduli spaces. It picks out $6$ dimensions and seems to be a natural point of entry into geometrical structures of interest to string theorists, some of whom have taken the point of view much further \cite{Witten},\cite{Gran}.

\subsection{Generalized $G_2$-structures}

A second case of open orbits, that of $ Spin(7,7)\times \R^*$ was considered by F. Witt in \cite{Wit}. We shall describe this briefly here since it has some parallels with the five-dimensional problem we will be mainly concerned with.

There is an open orbit on the $64$-dimensional half-spin representation which has stabilizer $G_2\times G_2$ where $G_2$ is the compact real form of the exceptional group. The invariant polynomial $f$ is of degree $8$ and somewhat difficult to write down \cite{Gyo}, but it is unnecessary to know it in detail. We realize the half-spinors as even forms.

The interpretation of a reduction of the structure group of $T\oplus T^*$, or our extension $E$, from $SO(7,7)$ to $G_2\times G_2$ can be seen in two steps. Since $G_2\times G_2\subset SO(7)\times SO(7)$, we have  a decomposition $E=V\oplus V^{\perp}$ into orthogonal subbundles where the induced inner product on $V$ is positive definite. This is already a generalized metric as defined above, and so comes equipped with two connections $\nabla^+,\nabla^-$ with torsion $\pm H$. The further reduction to $G_2\times G_2$ comes from the fact that $G_2$ is the stabilizer of a spinor. We obtain  spinor fields $\epsilon^+,\epsilon^-$, covariant constant with respect to $\nabla^+,\nabla^-$, equations already considered in the physics literature \cite{Gaunt}. 

As Witt shows, with respect to this decomposition of $E$, the form $\rho$ can be written as $\rho=e^f\epsilon^+\otimes \epsilon^-$ where $f$ is a function (the dilaton field). 

For our purposes it is interesting to note that $\hat\rho$ is obtained from a very natural procedure -- we define $R:E\rightarrow E$ to be the orthogonal reflection which is $1$ on $V$ and $-1$ on $V^{\perp}$. Its lift from $O(n,n)$ to $Pin(n,n)$ (it is orientation-reversing since $\dim V$ is odd) defines a map $\ast$ from even forms to odd forms and $\hat\rho=\ast\rho$. This map is closely related to the Hodge star operator and the equations $d\rho=d\!\ast\!\rho=0$ have an interpretation in terms of a Dirac equation satisfied by $\rho$. That, and the algebraic decomposability $\rho=e^f\epsilon^+\otimes \epsilon^-$, enable Witt to deduce that $\epsilon^\pm$ are covariant constant.
\vskip .25cm
The two examples here coming from $Spin(6,6)$ and $Spin(7,7)$ have very different properties: the generic case for $Spin(6,6)$ has a standard local normal form if we act via the group $\Diff(M) \ltimes \Omega_{closed}^2(M)$.  On the other hand $Spin(7,7)$ defines a metric whose curvature tensor gives local invariants. We shall investigate this aspect for the next case of $Spin(5,5)$

\section{A five-dimensional functional}
\subsection{Basic algebra}
 We consider now the geometrical implications of another group from Kimura and Sato's list --  $Spin(5,5)$. In fact $Spin(5,5)\times \R^*$ has an open orbit, but no invariant polynomial, and instead we consider the action of another group: $Spin(5,5)\times GL(2,\R)$ acting on two copies of the $16$-dimensional spin representation. We shall then derive a type of geometry on a five-manifold which belongs somehow to the same category of objects as those described above, but with distinctive features involving the Courant bracket.
\vskip .25cm
Let $S$ be the $16$-dimensional half-spin representation of $Spin(5,5)$. There are two of these but they are dual to each other. Consider  the standard action of the group  $Spin(5,5)\times GL(2,\R)$ on $S\otimes\R^2$. Then, according to \cite{Sat}, this has an open orbit with stabilizer $G_2^{(2)}\times SL(2,\R)$ where $G_2^{(2)}$ is the noncompact real form of the complex Lie group of type $G_2$. To see how this works, recall that $G_2^{(2)}\subset SO(3,4)$  and the adjoint representation gives a homomorphism $SL(2,\R)\rightarrow SO(2,1)$. We therefore have a homomorphism:
$$G_2^{(2)}\times SL(2,\R)\rightarrow  SO(3,4)\times SO(2,1)\subset SO(5,5)$$
which lifts to $Spin(5,5)$. We write $\R^{10}=V^{\perp}\oplus V$ corresponding  to this product of groups, with $V$ three-dimensional.

Under the action of $Spin(3,4)\times Spin(2,1)$, $S$ is expressed as a tensor product of spin representations of the factors: $\R^8\otimes \R^2$. Now  $G_2^{(2)}$  is the stabilizer of a non-null spinor $\psi$ in the $8$-dimensional spin representation  $Spin(3,4) \rightarrow SO(4,4)$ and $SL(2,\R)$ fixes a skew bilinear form $\epsilon \in \R^2\otimes\R^2$, thus $G_2^{(2)}\times SL(2,\R)$ stabilizes $$\psi\otimes \epsilon\in \R^8\otimes\R^2\otimes\R^2=S\otimes \R^2.$$
The dimension of $Spin(5,5)\times GL(2,\R)$ is $45+4=49$ and of $G_2^{(2)}\times SL(2,\R)$ $14+3=17$ and we have a $32$-dimensional open orbit.
In fact there are two orbits: the inner product on $V$  could have signature $(1,2)$ or $(2,1)$.
\vskip .25cm
The three-dimensional subspace $V$ of $\R^{10}$  on which $SO(2,1)$ acts generates a Clifford subalgebra which acts on $S$ through the $\R^2$ factor. Recall \cite{LM} that if $e_1,e_2,e_3$ is  a basis for $\R^3$ with $(e_1,e_1)=(e_2,e_2)=1, (e_3,e_3)=-1$ and $\omega=e_1\cdot e_2\cdot e_3$ then $\omega^2=1$ and commutes with everything: it decomposes the Clifford algebra into a sum of two copies of a real $2\times 2$ matrix algebra. In our situation this means that given a vector $v\in V\subset \R^{10}$, $v\cdot\omega:S\rightarrow S$ acts via a $2\times 2$ trace-free matrix $a(v)$ on the $\R^2$ factor of $S=\R^8\otimes\R^2$. The inner product $(v,v)$ obtained by restricting the form on $\R^{10}$ can be written as $\tr a(v)^2$.
\vskip .25cm
We describe next the invariant polynomial $f$ (see also \cite{Gyo}).
Clifford multiplication by $v\in \R^{10}$ maps the spinor space $S$ to the opposite spinor space, its dual $S^*$. We have the canonical pairing  $\langle \varphi,\psi\rangle$. Given $\varphi_1,\varphi_2\in S$ we can therefore define
$\langle v\cdot\varphi_1,\varphi_2\rangle\in \R$, which is  {\it symmetric} in $\varphi_1,\varphi_2$. Using the inner product,  define  the symmetric bilinear expression $P(\varphi_1,\varphi_2)\in \R^{10}$ by
$$(P(\varphi_1,\varphi_2),v)=\langle v\cdot\varphi_1,\varphi_2\rangle.$$

In particular we define $Q(\varphi)=P(\varphi,\varphi)$ (which is a null vector, indeed $S$ would have a quartic invariant $(Q(\varphi), Q(\varphi))$ otherwise).

Now consider the quartic real-valued function $f$ on $S\otimes \R^2$ defined by
\begin{equation}
f(\rho)=(Q(\rho_1),Q(\rho_2))
\label{deff}
\end{equation}
 where $\rho=(\rho_1,\rho_2)\in  S\otimes\R^2$. This is $Spin(5,5)$-invariant and a short calculation shows that if $A\in GL(2,\R)$ then
\begin{equation}
f(A(\rho))=(\det A)^2f(\rho)
\label{inv}
\end{equation}
and this is our relatively invariant polynomial.
\vskip .25cm
On the open orbit, we can see what $f$ is by considering $\rho=\psi\otimes \epsilon\in \R^8\otimes\R^2\otimes\R^2$, or, using a symplectic basis $u_1,u_2$ for $\R^2$, $\rho=\psi\otimes (u_1\otimes u_2-u_2\otimes u_1)$. Then 
$$\rho_1=-\psi\otimes u_2,\qquad \rho_2=\psi\otimes u_1.$$

Now the  spin representation of $Spin(3,4)$ has an invariant inner product and the spin representation of $Spin(2,1)=SL(2,\R)$ has the invariant skew form $\epsilon$. Thus the tensor product $S=\R^8\otimes\R^2$ has an invariant {\it skew}  form $B$. We use  Clifford multiplication by $\omega$ to identify $S$ and $S^*$, and then by invariance, the pairing of $S$ and $S^*$ must be a multiple of the skew form $B$ on $\R^8\otimes\R^2$.

This means that, up to a universal constant, 
$$\langle v\cdot \rho_1,\rho_1\rangle=B(v\cdot\omega\cdot(\psi\otimes u_2),\psi\otimes  u_2).$$
If $v\in V^{\perp}\subset S$ then $v$ acts only on the $8$-dimensional spin space -- the $\psi$ factor --  and since $\epsilon(u_2,u_2)=0$, this contributes zero to $\langle v\cdot \rho_1,\rho_1\rangle$. Thus it is only the $V$ component of a vector in $\R^{10}$ which contributes and then 
$$(Q(\rho_1),v)=\langle v\cdot \rho_1,\rho_1\rangle=B(v\cdot\omega\cdot (\psi\otimes  u_2),\psi\otimes u_2)=(\psi,\psi)\epsilon(a(v)u_2, u_2).$$
This 
gives  
$$Q(\rho_1)=(\psi,\psi)u_2\otimes u_2.$$
We identify the three-dimensional (adjoint) representation of $SL(2,\R)$ with   the symmetric part of $\R^2\otimes\R^2$. The symplectic basis $u_1,u_2$ generates the basis 
$$u_1\otimes u_1,\quad u_2\otimes u_2,\quad u_1\otimes u_2+u_2\otimes u_1$$
of the symmetric product which is the  basis of the Lie algebra ${\lie sl}(2,\R)$:
$$\pmatrix {0&1\cr
            0&0},\quad \pmatrix {0&0\cr -1& 0},\quad \pmatrix {-1&0\cr
            0&1}.$$
 Note that $Q(\rho_1)\in {\lie sl}(2,\R)$ is nilpotent and hence null, confirming our earlier statement.

These three basis elements correspond to $v_1,v_2,h$ where 
$$Q(\rho_1 +z\rho_2,\rho_1 +z\rho_2)=(\psi,\psi)\otimes (v_1+2zh+v_2z^2).$$
Now observe that
\begin{equation}
v_1\cdot\rho_1=0, \quad v_2\cdot\rho_2=0,\quad  v_1\cdot \rho_2=\omega\cdot\rho_1, \quad v_2\cdot\rho_1=-\omega\cdot\rho_2
\label{vrho}
\end{equation}

It follows from (\ref{deff})  that the quartic invariant is 
$$f(\rho)=(\psi,\psi)^2.$$

We need to compute the derivative of $f$ at $\rho$  --  this will define us the derivative of the volume form $\phi$ when we apply this algebra to the geometrical situation. The derivative will be an element $\rho'\in (S\otimes \R^2)^*$ such that $Df(\dot\rho)=\langle \rho',\dot\rho\rangle$. If $\rho$ lies in an open orbit, so does $\rho'$, and the roles of $S$ and $S^*$ are symmetrical in this respect.
\vskip .25cm
Now $f(\rho)=(Q(\rho_1),Q(\rho_2))$ and $Q(\varphi)=P(\varphi,\varphi)$ so differentiating we have 
$dQ(\dot\varphi)=2P(\varphi,\dot\varphi)$ and  
\begin{eqnarray*}
Df(\dot\rho)&=&2(P(\rho_1,\dot\rho_1),Q(\rho_2))+2(P(\rho_2,\dot\rho_2),Q(\rho_1))\\
&=&2[\langle Q(\rho_2)\cdot \rho_1,\dot\rho_1\rangle+\langle Q(\rho_1)\cdot \rho_2,\dot\rho_2\rangle]\\
&=& 2(\psi,\psi)[\langle \omega\cdot(\psi\otimes u_1),\dot\rho_1\rangle +\langle \omega\cdot(\psi\otimes u_2),\dot\rho_2\rangle]\\
&=& 2f^{1/2}[\langle \omega\cdot\rho_2,\dot\rho_1\rangle -\langle \omega\cdot\rho_1,\dot\rho_2\rangle]
\end{eqnarray*}
If we identify $\R^2$ with its dual using the skew form $\epsilon$ then this tells us that
\begin{equation}
\rho'=2f^{1/2}\omega\cdot \rho
\label{deriv}
\end{equation}
\vskip .5cm
To sum up, 
\begin{prp} \label{alg} Let $\rho\in S\otimes\R^2$ lie in an open orbit of $Spin(5,5)\times GL(2,\R)$, then
\begin{itemize}
\item
$\rho$ determines a $3$-dimensional subspace $V\subset \R^{10}$ on which the  quadratic form or its negative has signature $(2,1)$
\item
if $v\in V$ then the action of $v$  on $\rho$ by Clifford multiplication satisfies
$$v\cdot\rho=a(v)\omega\cdot\rho$$
where $a:V\rightarrow {\lie{sl}}(2,\R)$ is an isomorphism such that
$\tr a(v)^2=(v,v)$
\end{itemize}
\end{prp}

\subsection{The variational problem}
We now implant this algebra onto the differential geometry of a $5$-manifold $M$, following the ideas of \ref{open}. 

We replace the half-spin representations $S$ and $S^*$ of $Spin(5,5)$ by
$S=\Lambda^{ev}T^*$, $ S^*=\Lambda^{od}T^*$
and the  invariant Mukai pairing  is given by 
$$\langle \varphi,\psi\rangle=\varphi_0\wedge \psi_5-\varphi_2\wedge\psi_3+\varphi_4\wedge \psi_1\in \Lambda^5T^*.$$
The homogeneous quartic function $f$ can now be thought of as a map
$$f:\Lambda^{ev}T^*\rightarrow (\Lambda^5T^*)^2.$$
On the open orbit, $f$ is non-zero. Choosing an orientation on the $5$-manifold and taking the positive square root, we can define a volume form for each global section of $\Lambda^{ev}T^*\otimes \R^2$ which lies in the open orbit at each point:
$$\phi(\rho)=f^{1/2}(\rho)=\vert(\psi,\psi)\vert.$$
We shall  interpret the critical points of the variational problem using this volume form.
\vskip .5cm

Let $\rho\in \Omega^{ev}\otimes \R^2$ be an $\R^2$-valued  closed even form, lying in the open orbit at each point. We look for critical points of the volume functional 
$$V(\rho)=\int_M \!\!\phi(\rho)$$
on a given cohomology class in $H^{ev}(M,\R^2)$.
But $\phi=f^{1/2}$, so $D\phi=f^{-1/2}Df/2$ and from (\ref{deriv}) we see that $\rho'=\omega\cdot\rho$ and so the variational equation is 
\begin{equation}
d\rho=0\qquad d(\omega\cdot\rho)=0.
\label{max1}
\end{equation}

\begin{rmk} Note that $\rho\mapsto \omega\cdot \rho$ is the lift of the orthogonal reflection in $V^{\perp}$, so we are formally in a similar situation to generalized $G_2$ geometry.
Note further that since $\phi$ is homogeneous of degree $2$, Euler's formula gives
\begin{equation}
2\phi= D\phi(\rho)=\langle \omega\cdot\rho,\rho\rangle=2\langle \omega\cdot\rho_1,\rho_2\rangle.
\label{fidef}
\end{equation}
\end{rmk} 

\vskip .25cm
In our global setting we have a global $\R^2$-valued form $\rho$ and so  from Proposition \ref{alg}, the coefficients of $Q(\rho_1 +z\rho_2,\rho_1 +z\rho_2)=\phi\otimes (v_1+2zh+v_2z^2)$ define a three-dimensional space of global sections $v$ of $T\oplus T^*$ satisfying 
$$v\cdot\rho=a(v)\omega\cdot\rho.$$

This has a consequence for the interpretation of the variational equations.
\vskip .25cm
Since $a(v)$ is just a constant matrix, the equation $d(\omega\cdot\rho)=0$ in (\ref{max1}) tells us that
\begin{equation}
d(v\cdot\rho)= d(a(v)\omega\cdot\rho)=0.
\label{d1}
\end{equation}
Since $\omega^2=1$ and commutes with $v$ we also have for any two sections $v,w$,
$$(v\cdot w-w\cdot v)\cdot\rho=[a(v),a(w)]\rho$$
and since $d\rho=0$ it follows that 
\begin{equation}
d((v\cdot w-w\cdot v)\cdot\rho)=0.
\label{d2}
\end{equation}
Now, using (\ref{d1}) for $v$ and $w$, and (\ref{d2}), it follows from the definition of Courant bracket (\ref{courant}) that 
$$[v,w]\cdot\rho=0.$$
 This means that at each point $[v,w]$ lies in the annihilator of $\rho_1$ and $\rho_2$. But this is an {\it isotropic}  subspace preserved by $G_2^{(2)}\times SL(2,\R)$ and the $7$- and $3$-dimensional subspaces are the only invariant subspaces. Since the inner product is non-degenerate on these, we must have $[v,w]=0$. Hence: 
 
\begin{prp} A solution to the equations $d\rho=0=d\hat \rho$ gives three Courant-commuting sections of $T\oplus T^*$ whose inner products are constant.
\end{prp}

Picking a basis $v_i=X_i+\xi_i$ ($1\le i\le 3$), this means in particular that we have three commuting vector fields $X_1,X_2,X_3$. Moreover, the equation $d(v\cdot\rho)=0$ with $d\rho=0$ tells us that 
$$d(i_{X_i}\rho)+d\xi_i\wedge\rho=0$$
and  this can be written as 
\begin{equation}
{\mathcal L}_{X_i}\rho+d\xi_i\wedge\rho=0
\label{Lie}
\end{equation}
which says that $\rho$ (and similarly $\hat\rho=\omega\cdot\rho$) is preserved by a three-dimensional abelian subalgebra of the Lie algebra of $\Diff(M) \ltimes \Omega_{closed}^2(M)$.

\subsection{The geometrical structure}

We need to understand more concretely next how to obtain the three Courant-commuting sections of $T\oplus T^*$ from the pair of forms $\rho_1,\rho_2$. 
So for $A=1,2$, take $\rho_A=\varphi_0+\varphi_2+\varphi_4$, indexing the degrees of the forms by the subscript. We  use the  volume form $\phi$ and write $\varphi_4=i_{Y_A}\phi$ for a vector field $Y_A$. Since $d(i_{Y_A}\phi)=0$, this vector field is volume-preserving. We also write $\varphi_2=\omega$ and $\varphi_0=c$ to clarify their roles.

Then, for any $X+\xi$,
\begin{eqnarray*}
\langle (X+\xi)\cdot\rho_A,\rho_A \rangle &=& \langle i_{X}\omega+i_{X}i_{Y_A}\phi+c\xi+\xi\wedge\omega+\xi\wedge  i_{Y_A}\phi,c+ \omega+i_{Y_A}\phi\rangle\\
&=&i_{X}\omega \wedge i_{Y_A}\phi-(i_{X}i_{Y_A}\phi+\xi\wedge\omega)\wedge\omega+2c\xi\wedge i_{Y_A}\phi
\end{eqnarray*}
Now $i_X\omega\wedge \phi=0$ since it is a $6$-form in $5$ dimensions. This means that 
$$i_{X}\omega\wedge i_{Y_A}\phi=-(i_{X}i_{Y_A}\omega)\phi$$
Thus we can rewrite the above as
\begin{equation}
\langle (X+\xi)\cdot\rho_A,\rho_A \rangle= -2i_Xi_{Y_A}\omega \phi-\xi\wedge\omega^2 +2c\xi\wedge i_{Y_A}\phi.
\label{one}
\end{equation}
By definition this is $(Q(\rho_A),X+\xi)$ so with  $Q(\rho_A)=\phi\otimes(X_A+\xi_A)=\phi\otimes v_A$ we have
$$(Q(\rho_A),X+\xi)=(i_{X_A}\xi+i_X\xi_A)\phi/2.$$
Since also $\xi\wedge \phi=0$, we have
$\xi\wedge i_{X_A}\phi=(i_{X_A}\xi) \phi$ and then we rewrite this as
\begin{equation}
(Q(\rho_A),X+\xi)=\xi\wedge i_{X_A}\phi/2+i_X\xi_A\phi/2
\label{two}
\end{equation}
Comparing (\ref{one}) and (\ref{two}), we 
 get
\begin{equation}
\xi_A= -4i_{Y_A}\omega
\label{x1}
\end{equation}
 and
\begin{equation}
i_{X_A}\phi=-2\omega^2 +4ci_{Y_A}\phi
\label{X1}
\end{equation}

\begin{rmk} \label{volpreserve} If only $d\rho_A=0$ holds, then $c$ is a constant, $\omega$ and $i_{Y_A}\phi$ are closed, and hence from (\ref{X1}) $d(i_{X_A}\phi)=0$ so that $X_A$ is volume-preserving. 
\end{rmk}

\vskip .5cm

We try next to use the transformations in $\Diff(M) \ltimes \Omega_{closed}^2(M)$ to reduce the data to a more manageable form. If $\rho_1=\varphi_0+\varphi_2+\varphi_4$ and $\rho_2=\psi_0+\psi_2+\psi_4$ in terms of even forms, then the degree zero part is closed and hence constant so by an $SL(2,\R)$ transformation we can assume that $\varphi_0=0$.

 Assume that we are in the generic case $\psi_0\ne 0$, then by another transformation we can take it to be $1$. But then transforming by $\exp(-\psi_2)$ ($\psi_2$ is closed) we can make $\psi_2=0$. Thus after the action of $SL(2,\R)$  and a closed $B$-field we can take
\begin{equation}
\rho_1=\omega+i_{Y_1}\phi,\quad \rho_2=1+i_{Y_2}\phi
\label{oo}
\end{equation}
for vector fields $Y_1,Y_2$. It then follows from (\ref{X1}) that $$i_{X_1}\phi=-2\omega^2, \quad i_{X_2}\phi= 4i_{Y_2}\phi$$
 so $X_2=4Y_2$, and from (\ref{x1})
$$\xi_1= -4i_{Y_1}\omega,\quad \xi_2=0$$
Since 
$\phi^2=(Q(\rho_1), Q(\rho_2))=(i_{X_1}\xi_2+i_{X_2}\xi_1)\phi^2/2$, we then see that $$2=i_{X_2}\xi_1=-4i_{X_2}i_{Y_1}\omega=-16\omega(Y_1,Y_2).$$
\vskip .25cm
The integrability condition $d\hat\rho=0$ is 
\begin{equation}
d(\omega\cdot \rho_1)=0=d(\omega\cdot \rho_2)
\label{int}
\end{equation}
 From (\ref{vrho})  
 $v_1\cdot \rho_2=\omega\cdot \rho_1$ so 
$d((X_1+\xi_1)\cdot \rho_2)=0$ and similarly $d((X_2+\xi_2)\cdot \rho_1)=0.$
Since $\xi_2=0$, the second equation in (\ref{int}) gives
$${\mathcal L}_{X_2}\omega=0,\quad {\mathcal L}_{X_2}(i_{Y_1}\phi)=0.$$
Since $X_2=4Y_2$ and 
${\mathcal L}_{X_2}(i_{Y_1}\phi)=i_{[X_2,Y_1]}\phi$ this is equivalent to 
$${\mathcal L}_{Y_2}\omega=0,\quad [Y_1,Y_2]=0.$$

From the first equation in (\ref{int}) we have 
${\mathcal L}_{X_1}(i_{Y_2}\phi)=0$ which is equivalent to $[X_1,Y_2]=0$. But $i_{X_1}\phi=-2\omega^2$ so this is already implied by ${\mathcal L}_{Y_2}\omega=0$. The new relation is   $d\xi_1=0$ which from (\ref{x1}) is equivalent to  ${\mathcal L}_{Y_1}\omega=0$.

Therefore we  obtain the following:

\begin{thm} \label {generic} Let $\rho$ be a solution of the  variational equations $d\rho=0=d\hat\rho$ such that the degree zero component of $\rho$ is non-zero. Then up to  the action of $SL(2,\R)$ and a closed $2$-form $B$, the structure is equivalent to:
\begin{itemize}
\item
a closed $2$-form $\omega$ and a volume form $\phi$
\item
two commuting vector fields $Y_1,Y_2$, preserving  $\omega$ and $\phi$, such that  
\item
$\omega(Y_1,Y_2)=-1/8.$ 
\end{itemize}
The form $\rho$ then has the invariant volume form $\phi$, and is defined by 
$$\rho_1=\omega+i_{Y_1}\phi,\quad \rho_2=1+i_{Y_2}\phi.$$
\end{thm}
\vskip .5cm

\begin{rmk} By considering $Q(\rho_1+\rho_2)$ we can easily find the third of the Courant-commuting sections of $T\oplus T^*$. We have altogether
$$Y_1-i_{Y_2}\omega,\quad Y_2,\quad X-i_{Y_1}\omega$$
where $X$ is defined by $i_X\phi=-8\omega^2$.
\end{rmk}

\begin{ex} Clearly from Theorem \ref{generic} we can find compact examples e.g. $S^3\times T^2$ where $T^2$ is a flat symplectic torus and $S^3$ has a fixed  volume form. In this case $Y_1,Y_2$ are the linear vector fields on the torus and $X=0$. Or we could take $S^1\times T^2\times S^2$ where $\omega$ is a symplectic form on the four-manifold $T^2\times S^2$ -- here $X$ is a vector field on the circle factor.
\end{ex}

\subsection{A normal form}\label{normform}

Under a further regularity condition, Theorem \ref{generic}  allows us to give a local normal form: 
\begin{prp} If the $2$-form $\omega$ has maximal rank, then locally $\rho$ is equivalent to 
\begin{itemize}
\item
$\rho_1=dx_1\wedge dx_2+dx_3\wedge dx_4+dx_1\wedge dx_3\wedge dx_4\wedge dx_5$
\item
$\rho_2= 1+dx_2\wedge dx_3\wedge dx_4\wedge dx_5$
\end{itemize}
\end{prp}

\begin{prf} We use  Theorem \ref{generic}.  Since ${\mathcal L}_{Y_1}\omega=0={\mathcal L}_{Y_2}\omega$ and $\omega$ is closed, the $1$-forms $i_{Y_1}\omega,  i_{Y_2}\omega$ are closed, and since $\omega(Y_1,Y_2)\ne 0$, they are everywhere linearly independent.  So locally $i_{Y_A}\omega =dx_A$ and $(x_1,x_2)$ defines a map $p:M\rightarrow \R^2$ with fibre a $3$-manifold $N$.
Since $\omega(Y_1,Y_2)=-1/8$,  then 
$$-8dp(aY_1+bY_2)=(b,-a)$$
 so $Y_1,Y_2$ cover the standard commuting vector fields on $\R^2$. But $Y_1,Y_2$ are themselves  commuting vector fields on $M$, so they define a flat $\Diff(N)$ connection over the open set in $\R^2$. This can be trivialized by a local diffeomorphism and then we have coordinates $x_1,\dots, x_5$ with
$Y_1=c_1\partial/\partial x_2, 
Y_2= c_2 \partial/\partial x_1.$  

Define $\omega_0=\omega+i_{Y_1}\omega\wedge i_{Y_2}\omega$, then $i_{Y_1}\omega_0=i_{Y_2}\omega_0=0$. Thus $\omega$ is determined by $\omega_0$ -- its restriction to a fibre. It is also closed and invariant by $Y_1,Y_2$, thus $\omega = \omega_0+cdx_1\wedge dx_2$ where $\omega_0$ is  a closed $2$-form on an open set $U\subset \R^3$. Rescale the coordinates $x_1,x_2$ so that $c=1$.

 If $\omega$ has rank $4$, then $\omega_0$ has rank $2$. Now the volume form $\phi$ together with $dx_1\wedge dx_2$ defines a volume form $\nu$ on the $3$-dimensional fibres. Moreover, $Y_1,Y_2$ preserve $\phi$ and $dx_1\wedge dx_2$, and so in the product decomposition of $M$, $\nu$ is independent of $x_1,x_2$.
 
We  write $\omega_0=i_X\nu$ for a non-vanishing volume-preserving vector field $X$ on $U$.  Since $i_X\omega=0$ and $\omega_0$ is closed, ${\mathcal L}_X\omega_0=0$, and so $\omega_0$ is pulled back from the quotient of $U$ by the action generated by $X$. We can therefore write $\omega_0=dx_3\wedge dx_4$ and $X=\partial/\partial x_5$. 
The volume form $\nu$ on $U$ must now be  $dx_5\wedge dx_3\wedge dx_4$.

Hence we have a normal form 
\begin{eqnarray*}
\rho_1&=&dx_1\wedge dx_2+dx_3\wedge dx_4+dx_1\wedge dx_3\wedge dx_4\wedge dx_5\\
\rho_2&=& 1+dx_2\wedge dx_3\wedge dx_4\wedge dx_5
\end{eqnarray*}
and the space of sections $V$ is spanned by 
$$\left\{\frac{\partial}{\partial x_5}+dx_1,\quad \frac{\partial}{\partial x_1},\quad \frac{\partial}{\partial x_2}+dx_2\right\}$$
\end{prf}

\section{The six-dimensional geometry}

\subsection{A Hamiltonian flow}

 The affine space of closed forms in a fixed cohomology class in $H^{ev}(M,\R)$ has a natural constant metric: if $\dot\varphi=d\alpha$ then 
$$(\dot\varphi,\dot\varphi)=\int_M \!\langle d\alpha,\alpha\rangle.$$
Using the skew form $\epsilon$ on $\R^2$ we can therefore give the forms in a cohomology class $a\in H^{ev}(M,\R^2)$ a symplectic structure $\Omega$. We used $\epsilon$ before to define $\hat\rho\in \Omega^{od}\otimes \R^2$, and using the same notation we write 
$$\Omega(\dot\rho,\dot\sigma)=\int_M\! \!\langle \dot\rho,\beta\rangle$$
where $\dot\sigma=d\beta$.

The volume functional $V(\rho)$ generates a Hamiltonian flow which is given by the equation 
$$\Omega(\dot\rho,d\beta)=DV(d\beta)=\int_M\!\!\langle \hat\rho, d\beta\rangle.$$
But by definition 
$$\Omega(\dot\rho,d\beta)=\int_M \!\!\langle \dot\rho,\beta\rangle$$
so we obtain 
\begin{equation}
\frac{\partial \rho}{\partial t}=d\hat\rho
\label{evolve}
\end{equation}
or equivalently
$$\frac{\partial \rho_1}{\partial t}=d(\omega\cdot\rho_1),\quad \frac{\partial \rho_2}{\partial t}=d(\omega\cdot\rho_2).$$
\vskip .25cm
Consider as earlier, for an indeterminate $z$, $\rho(z)=\rho_1+z\rho_2$. Then
$$\omega\cdot\rho(z)=v_1\cdot\rho_2-zv_2\cdot \rho_1=(z^{-1}v_1-zv_2)\cdot \rho(z)$$
since $v_A\cdot\rho_A=0$. We set $u(z)=z^{-1}v_1-zv_2$.

Now
$$Q(\rho(z),\rho(z))=\phi\otimes (v_1+2zh+v_2z^2)=\phi\otimes v(z)$$ and $v(z)\cdot\rho(z)=0$. 

Note that because $v(z)$ annihilates $\rho(z)$ we could replace $u(z)$ by $u(z)-z^{-1}v(z)$ to get a polynomial in $z$ and by $u(z)+z^{-1}v(z)$ to get a polynomial in $z^{-1}$, so what follows is also valid for $z=0$ and $\infty$.

Differentiating $v(z)\cdot\rho(z)=0$, we find
$$\frac{\partial v(z)}{\partial t}\cdot\rho(z)=-v(z)\cdot d(\omega\cdot \rho(z))=-v(z)\cdot d(u(z)\cdot \rho(z)).$$
  We have  $v(z)\cdot\rho(z)=0$ and since $\rho(z)$ is closed $d((u(z)\cdot v(z)-v(z)\cdot u(z))\cdot \rho(z))=0$ (see (\ref{d2})). Now use the Courant bracket formula (\ref{courant}) to obtain
 $$\left(\frac{\partial v(z)}{\partial t}+[v(z),u(z)]\right)\cdot\rho(z)=0.$$
By $G_2^{(2)}$ invariance this means that 
\begin{equation}
\frac{\partial v(z)}{\partial t}+[v(z),u(z)]=\lambda v(z)
\label{evol}
\end{equation}
and since the left hand side is quadratic in $z$, as is $v(z)$, $\lambda$ is independent of $z$, and so equating coefficients we obtain:
\begin{eqnarray*}
v_1'&=&-2[h,v_1]+\lambda v_1\\
h'&=&[v_1,v_2]+\lambda h\\
v_2'&=&2[h,v_2]+\lambda v_2
\end{eqnarray*}
(Note that since $(h,h)=2$, we have
$2\lambda=-(h,[v_1,v_2]).$)

These equations can be written as 
\begin{equation}
v_i'=c_{ijk}[v_j,v_k]+\lambda v_i
\label{nm}
\end{equation}
where $c_{ijk}$ are the structure constants of ${\lie sl}(2,\R)$. 
\vskip .25cm
If we used the structure constants $\epsilon_{ijk}$ of ${\lie su}(2)$, and took $\lambda=0$, with $v_i(t)$  taking values in a finite-dimensional Lie algebra, then we would have {\it Nahm's equations}. In fact, in that case, by rescaling and reparametrizing $t$, we can  transform (\ref{nm}) with an arbitrary function $\lambda(t)$ into Nahm's form.

When the $v_i$ are vector fields on a {\it three}-manifold $N$, with the bracket the Lie bracket, and $\lambda=0$, Nahm's equations describe a hypercomplex structure on the four-manifold $\R\times N$ \cite{Hit5}, which generalizes the original observation of Ashtekar et al \cite{Ash} that if the $v_i$ are volume-preserving, one obtains a hyperk\"ahler manifold. Clearly the change from $\epsilon_{ijk}$ to $c_{ijk}$ alters these constructions to yield the signature $(2,2)$ versions of hypercomplex and hyperk\"ahler. In the case where  $\lambda\ne 0$ and is a function of both $t$ and the three-manifold $N$, we cannot absorb $\lambda$, but on the other hand Nahm's equations correspond to the situation where time $t$ is a harmonic function on $\R\times N$ (see \cite{Hit5}). The evolution equations in the more general case give the non-zero $\lambda$  term.
\vskip .25cm
We see from these remarks that the shape of (\ref{nm}) is by no means unfamiliar. We have the Courant bracket, but the vector field part of this is just the Lie bracket, and this would give us a signature $(2,2)$ solution to Einstein's equations if the vector fields were tangent to a $3$-dimensional foliation. This is not the case in general, but as we shall see next, the six-dimensional geometrical structure being generated has something in common with four-dimensional self-duality.
\subsection{The six-dimensional structure}

On the six-manifold $\R\times M$ consider the following two sections of $T\oplus T^*$:
$$v(z),\quad w(z)=\frac{\partial}{\partial t}-2dt-u(z).$$
\begin{prp} The rank $2$ subbundle $E_z\subset T\oplus T^*$ spanned at each point by $v(z),w(z)$  is isotropic and Courant-integrable.
\end{prp}
(By this we mean that the Courant bracket of any two sections of $E_z$ lies in $E_z$ -- because of the property (\ref{c1}) of the Courant bracket, this only makes sense if $E_z$ is isotropic.)
\begin{prf}
From the definition of $v(z)$, we have $(v(z),v(z))=0$. Also  
$$(u(z),v(z))=(z^{-1}v_1-zv_2, v_1+2zh+v_2z^2)=z(v_1,v_2)-z(v_2,v_1)=0.$$
We have
$(u(z),u(z))=-2(v_1,v_2)=2$, using $(v_1,v_2)=\tr a(v_1)a(v_2)$. But 
$$\left(\frac{\partial}{\partial t}-2dt,\frac{\partial}{\partial t}-2dt\right)=-2,$$ and so $w(z)$ is null.
Hence all inner products of $v(z),w(z)$ vanish.
\vskip .25cm
From (\ref{evol}) we get 
$$\left[\frac{\partial}{\partial t}-2dt+u(z), v(z)\right]=\lambda v(z)$$
so for each $z$ we get a $2$-dimensional Courant-integrable subbundle $E_z\subset T\oplus T^*$. 
\end{prf}

Now since $(u(z),u(z))=2$ the vector field part of $u(z)$ is always non-zero, and the vector field part of $w(z)$ has a $\partial/\partial t$ component. Thus for each $z$ we get a foliation of $\R\times M$ by surfaces, a $5$-dimensional family in all. These are analogues of the $3$-dimensional family of $\alpha$-planes in self-dual $4$-manifolds, and indeed there must be a description of this structure via a $5$-dimensional twistor space, though we shall not investigate that here.  
Instead we look for another characterization of this structure, and for this we compare with the evolution equation for an ordinary $G_2$ structure.
\vskip .25cm
In \cite{Hit2}, the open orbit of $GL(7,\R)$ on $3$-forms or $4$-forms was used to characterize metrics with holonomy $G_2$ in a similar fashion (but in the non-generalized setting) to the above. In that case a gradient flow equation for $4$-forms in a  fixed cohomology class was shown to generate a metric with $Spin(7)$ holonomy on $\R\times M^7$. In a similar way, Witt in \cite{Wit} showed how this could be done in a generalized setting, obtaining $Spin(7)\times Spin(7)$ structures -- connections $\nabla^+,\nabla^-$ with skew torsion in  $8$ dimensions with covariant constant spinors. The presence of the group $G_2^{(2)}$ in our case suggests that our six-manifold should have a structure related to $Spin(3,4)$. This is not in itself an open orbit question -- a manifold of holonomy $Spin(7)$ is determined by a differential form of degree $4$ which is not stable -- it has a special algebraic form. On the other hand, the  ``integrability" condition is just the single equation that  it should be closed. This is our approach here.
\vskip .5cm
The equation
$$\frac{\partial \rho}{\partial t}=d\hat\rho$$
together with $d\rho=0$ is equivalent to the condition that 
$$\sigma(z)=dt\wedge \hat\rho(z)+\rho(z)$$
should be a closed form on $\R\times M$. Here $\sigma(z)=\sigma_1+z\sigma_2$ where
$$\sigma_1=dt\wedge v_1\cdot \rho_2+\rho_1,\quad \sigma_2=-dt\wedge v_2\cdot \rho_1+\rho_2.$$
We need to understand the algebraic properties of the pair of forms $(\sigma_1,\sigma_2)$.
\begin{prp} $v(z)\cdot\sigma(z)=w(z)\cdot\sigma(z)=0.$ 
\end{prp}
\begin{prf}
We can write $\sigma(z)=dt\cdot u(z)\cdot \rho(z)+\rho(z)$, and then
$$v(z)\cdot\sigma(z)=-dt\cdot v(z)\cdot u(z)\cdot \rho(z)+v(z)\cdot\rho(z)=dt\cdot u(z)\cdot v(z)\cdot \rho(z)+v(z)\cdot\rho(z)=0$$
using the orthogonality of $u(z)$ and $v(z)$ and $v(z)\cdot \rho(z)=0$.
Now
$$u(z)\cdot\sigma(z)=-dt\cdot u(z)^2\cdot \rho(z)+u(z)\cdot \rho(z)=-2dt\cdot\rho(z)+u(z)\cdot\rho(z)$$
since $(u(z),u(z))=2$,
and 
$$\left(\frac{\partial}{\partial t}-2dt\right)\cdot\sigma(z)=u(z)\cdot \rho(z)-2dt\cdot \rho(z)$$
and putting these together gives $w(z)\cdot\rho(z)=0$.
\end{prf}

This proposition tells us that each form $\sigma(z)=\sigma_1+z\sigma_2$ has a two-dimensional annihilator (in fact no more than two, as can be seen by using the normal forms in \ref{normform}). More than that, we see that, as $z$ varies, the two-dimensional subspaces span a $4$-dimensional space with basis
$$\frac{\partial}{\partial t}-2dt,\quad v_1, \quad v_2,\quad h.$$
The inner products of these basis vectors form a constant matrix of signature $(2,2)$. Thus on our six-manifold, $T\oplus T^*$ has a distinguished rank $4$ trivial subbundle $V$, with nondegenerate inner product, so it has an orthogonal complement and the structure group reduces to $SO(4,4)\times SO(2,2)$. Let $S^+, S^-$ be the two half-spin representations of $Spin(4,4)$ and $\Delta^+,\Delta^-$ those of $Spin(2,2)=SL(2,\R)\times SL(2,\R)$. Then $\sigma_1,\sigma_2$ lie in
$$S^+\otimes \Delta^+ + S^-\otimes \Delta^-.$$
But $\sigma_1+z\sigma_2$ is annihilated by an isotropic $2$-dimensional subspace of this four-dimensional space, which puts the vector space spanned by $\sigma_1,\sigma_2$ in correspondence with one of the spin spaces of $Spin(2,2)$, say $\Delta^+$. The pair can therefore be considered as an element of $S^+\otimes \Delta^+\otimes \R^2$ and we are looking at its stabilizer. The stabilizer of a non-null vector $\psi\in S^+$ is $Spin(3,4)$, using triality, so we have  the pair $(\sigma_1,\sigma_2)$ given by
$$\psi\otimes \epsilon \in S^+\otimes \Delta^+\otimes \Delta^+.$$
In other words, we can describe the geometry on the six-manifold as a closed $\R^2$-valued differential form whose stabilizer in $Spin(6,6)\times GL(2,\R)$ is the subgroup $Spin(3,4)\times SL(2,\R)\times SL(2,\R)$.

\begin{rmk} There are really two types of structure depending on whether we take the forms $\sigma_A$ to be even or odd. The same situation holds   for $Spin(7)\times Spin(7)$ structures in \cite{Wit} and \cite{Gaunt}.
\end{rmk}
 
\vskip 1cm
 Mathematical Institute, 24-29 St Giles, Oxford OX1 3LB, UK
 
 hitchin@maths.ox.ac.uk
\end{document}